\documentclass{article}
\usepackage[english,russian,ukrainian]{babel}
\usepackage{latexsym,amssymb}
\usepackage[cp1251]{inputenc}
\usepackage{amsfonts,amssymb}
\usepackage{euscript}
\usepackage{latexsym,amssymb,amsmath}
\usepackage[cp1251]{inputenc}
\usepackage{amsfonts,amssymb}
\usepackage{graphicx,graphics,hhline}
\usepackage{euscript}
\newcounter{theorem} 
\newcounter{lemma} 
\renewcommand{\thetheorem}{\arabic{theorem}}
\renewcommand{\thelemma}{\arabic{lemma}}
\newcommand{\theor}{\par\refstepcounter{theorem}%
{\bf Теорема \thetheorem .}\,\,}
\newcommand{\lem}{\par\refstepcounter{lemma}%
{\bf Лемма \thelemma .}\,\,}

\,
\,

\begin{document}

\vspace*{7mm}

\LARGE

\addtocounter{page}{-1}
\thispagestyle{empty}

\begin{center}
\textbf{Monogenic functions in
finite-dimensional\\
 commutative associative algebras}
\end{center}
\vskip4mm
\large
\begin{center}\textbf{V.\,S.\,Shpakivskyi}\end{center}

\vspace{7mm}
\small
Let $\mathbb{A}_n^m$ be an arbitrary $n$-dimensional
commutative associative algebra
over the field of complex numbers with $m$ idempotents. Let
$e_1=1,e_2,\ldots,e_k$ with $2\leq k\leq 2n$ be elements of $\mathbb{A}_n^m$ which are linearly
independent over the field of real numbers. We consider monogenic
(i.~e. continuous and differentiable in the sense of Gateaux)
functions of the variable $\sum_{j=1}^k x_j\,e_j$\,, where $x_1,x_2,\ldots,x_k$ are
real, and obtain a constructive description of all mentioned
functions by means of holomorphic functions of complex variables.
It follows from this description that monogenic functions have
Gateaux derivatives of all orders. The present article is generalized of the author's paper \cite{Shpakivskyi-2014}, where mentioned results are obtained for $k=3$.
\vspace{7mm}

\large

\section{Introduction}

Apparently, W.~Hamilton (1843) made the first attempts to
construct an algebra associated with the three-dimensional Laplace
equation\medskip
\begin{equation} \label{Lap3}
 \Delta_3 u(x,y,z):=
\left(\frac{{\partial}^2}{{\partial x}^2}+
 \frac{{\partial}^2}{{\partial y}^2}+
\frac{{\partial}^2}{{\partial z}^2}\right)u(x,y,z)
 =0\,
\end{equation}\medskip
meaning that components of hypercomplex functions satisfy the
equation
(\ref{Lap3}).
He constructed an algebra of noncommutative quaternions over the
field of real numbers $\mathbb{R}$ and made a base for developing
the hypercomplex analysis.

C.~Segre \cite{Segre} constructed an algebra of commutative
quaternions over the field $\mathbb{R}$ that can be considered as
a two-dimensional commutative semi-simple algebra of bicomplex
numbers over the field of complex numbers $\mathbb{C}$.
M.~Futagawa \cite{Futugawa} and J.~Riley \cite{Riley} obtained a constructive
description of analytic function of a bicomplex variable, namely,
they proved that such an analytic function can be constructed with
an use of two holomorphic functions of complex variables.

 F.~Ringleb \cite{Ringleb} and S.~N.~Volovel'skaya \cite{Volovel'skaya-1939, Volovel'skaya-1948} succeeded in developing a function
theory for noncommutative algebras with unit over the
real or complex fields, by pursuing a definition of the differential of
a function on such an algebra suggested by Hausdorff in \cite{Hausdorff}. These definitions
make the \textit{a priori} severe requirement that the coordinates of the
function have continuous first derivatives with respect to the coordinates
of the argument element. Namely, F.~Ringleb \cite{Ringleb} considered an arbitrary
finite-dimensional associative (commutative or not)  \textit{semi-simple}
algebra over the field
$\mathbb{R}$. For given class of functions which maps the mentioned
algebra onto itself, he obtained a constructive description by
means of real and complex analytic functions.

S.~N.~Volovel'skaya developed the Hausdorff's idea defining the \textit{monogenic} functions on
 \textit{non-semisimple} associative algebras and she generalized the Ringleb's results for such algebras.
    In the paper \cite{Volovel'skaya-1939} was obtained a
constructive description of monogenic functions in a special three-dimensional non-commutative
 algebra over the field $\mathbb{R}$. The results of paper \cite{Volovel'skaya-1939} were generalized in the paper \cite{Volovel'skaya-1948} where Volovel'skaya obtained a
constructive description of monogenic functions in non-semisimple associative algebras of the first category over $\mathbb{R}$.

A relation between spatial potential fields and analytic functions
given in commutative algebras was established by P.~W.~Ketchum
\cite{Ketchum-28} who shown that every analytic function
$\Phi(\zeta)$ of the variable $\zeta=xe_1+ye_2+ze_3$ satisfies the
equation (\ref{Lap3}) in the case where the elements $e_1, e_2,
e_3$ of a commutative algebra satisfy the condition
\begin{equation}\label{garmonichnyj_bazys-ogljad}
     e_1^2+e_2^2+e_3^2=0\,,
\end{equation}
because
\begin{equation}\label{garm}
\frac{{\partial}^{2}\Phi}{{\partial x}^{2}}+
\frac{{\partial}^{2}\Phi}{{\partial y}^{2}}+
\frac{{\partial}^{2}\Phi}{{\partial z}^{2}}\equiv{\Phi}''(\zeta) \
(e_1^2+e_2^2+e_3^2)=0\,,\medskip
\end{equation}
where $\Phi'':=(\Phi')'$ and $\Phi'(\zeta)$ is defined by the
equality $d\Phi=\Phi'(\zeta)d\zeta$.

We say that a commutative associative algebra $\mathbb A$ is {\it
harmonic\/} (cf.
\cite{Ketchum-28,Mel'nichenko75,Plaksa}) if in
$\mathbb A$ there exists a triad of linearly independent vectors
$\{e_1,e_2,e_3\}$ satisfying the equality
(\ref{garmonichnyj_bazys-ogljad}) with $e_k^2\ne 0$ for
$k=1,2,3$. We say also that such a triad $\{e_1,e_2,e_3\}$ is {\it
harmonic}.

P.~W.~Ketchum \cite{Ketchum-28} considered the C.~Segre algebra
of quaternions \cite{Segre} as an example of harmonic algebra.

Further M.~N.~Ro\c{s}cule\c{t} establishes a relation between monogenic functions in commutative algebras and partial differential equations. He defined \textit{monogenic} functions $f$ of the variable $w$ by the equality $df(w)\,dw=0$.
 So, in the paper \cite{Rosculet} M.~N.~Ro\c{s}cule\c{t} proposed a procedure for constructing an infinite-dimensional topological vector space with commutative multiplication such that monogenic functions in it are the \textit{all} solutions of the equation
\begin{equation}\label{intr}
\sum\limits_{\alpha_0+\alpha_1+\ldots+\alpha_p=N}C_{\alpha_0,\alpha_1,\ldots,\alpha_p}\,
\frac{\partial^N \Phi}{\partial x_0^{\alpha_0}\,\partial x_1^{\alpha_1}\,\ldots\partial x_p^{\alpha_p}}=0,
\end{equation}
with $C_{\alpha_0,\alpha_1,\ldots,\alpha_p}\in\mathbb{R}$.
In particular, such infinite-dimensional topological vector space are constructed for the Laplace
equation (\ref{garm}).  In the paper \cite{Rosculet-56} Ro\c{s}cule\c{t} finds a certain connection between monogenic functions in commutative algebras and systems of partial differential equations.

I.~P.~Mel'nichenko  proposed for describing solutions of the equation (\ref{intr})
to use hypercomplex functions differentiable in the sense of Gateaux, since in this case the conditions of monogenic are the least restrictive. He started to implement this approach with respect to the thee-dimensional Laplace equation (\ref{garm}) (see \cite{Mel'nichenko75}).
 Mel'nichenko proved that there exist exactly $3$
three-dimensional harmonic algebras with unit over the field
$\mathbb C$ (see
\cite{Mel'nichenko75,Melnichenko03,Plaksa}).

In the paper \cite{Pogorui-Rod-Shap}, the authors develop the Melnichenko's idea  for the equation (\ref{intr}), and considered several examples.

The investigation of partial differential equations using the hypercomplex methods is effective if hypercomplex monogenic (in any sense) functions can be constructed explicitly. On this way the following results are obtained.

Constructive descriptions of monogenic (i.~e. continuous and
differentiable in the sense of Gateaux) functions taking values in
the mentioned three-dimensional harmonic algebras by means three
corresponding holomorphic functions of the complex variable are
obtained in the papers \cite{Pl-Shp1,Pl-Pukh,Pukh}. Such
descriptions make it possible to prove the infinite
differentiability in the sense of Gateaux of monogenic functions
and integral theorems for these functions that are analogous to
classical theorems of the complex analysis (see, e.~g.,
\cite{Pl-Shp3,Plaksa12}).

Furthermore, constructive descriptions of monogenic functions
taking values in special $n$-dimensional
commutative algebras by means $n$ holomorphic functions of complex
variables are obtained in the papers \cite{Pl-Shp-Algeria,
Pl-Pukh-Analele}.

In the paper \cite{Shpakivskyi-2014}, by author is obtained a constructive description of \textit{all} monogenic
functions of the variable $x_1e_1+x_2e_2+x_3e_3$ taking values in an arbitrary $n$-dimensional
commutative associative algebra with unit by means of holomorphic
functions of complex variables.
It follows from this description
that monogenic functions have Gateaux derivatives of all orders.

In this paper we extend the results of the paper \cite{Shpakivskyi-2014} to
 monogenic functions of the variable $\sum\limits_{r=1}^k x_re_r$, where $2\leq k \leq 2n$.

\section{The algebra $\mathbb{A}_n^m$}

Let $\mathbb{N}$ be the set of natural numbers.
We fix the numbers $m,n\in\mathbb{N}$ such that $m\leq n$.
Let $\mathbb{A}_n^m$ be an arbitrary commutative associative algebra with
 unit over the field of complex number $\mathbb{C}$.
   E.~Cartan \cite[p.~33]{Cartan}
  proved that there exist a basis $\{I_r\}_{r=1}^{n}$ in $\mathbb{A}_n^m$
  satisfying the following multiplication rules:
  \vskip3mm
1.  \,\,  $\forall\,r,s\in[1,m]\cap\mathbb{N}\,:$ \qquad
$I_rI_s=\left\{
\begin{array}{rcl}
0 &\mbox{if} & r\neq s,\vspace*{2mm} \\
I_r &\mbox{if} & r=s;\\
\end{array}
\right.$

\vskip5mm

2. \,\,  $\forall\,r,s\in[m+1,n]\cap\mathbb{N}$\,: \qquad $I_rI_s=
\sum\limits_{p=\max\{r,s\}+1}^n\Upsilon_{r,p}^{s}I_p$\,;

\vskip5mm

3.\,\, $\forall\,s\in[m+1,n]\cap\mathbb{N}$\,\, 
 $\exists!\;
 u_s\in[1,m]\cap\mathbb{N}$ \,$\forall\,
r\in[1,m]\cap\mathbb{N}$\,:\,\,

\begin{equation}\label{mult_rule_3}
I_rI_s=\left\{
\begin{array}{ccl}
0 \;\;\mbox{if}\;\;  r\neq u_s\,,\vspace*{2mm}\\
I_s\;\;\mbox{if}\;\;  r= u_s\,. \\
\end{array}
\right.\medskip
\end{equation}
Moreover, the structure constants
$\Upsilon_{r,p}^{s}\in\mathbb{C}$
 satisfy the associativity conditions:
\vskip2mm (A\,1).\,\, $(I_rI_s)I_p=I_r(I_sI_p)$ \,\,
$\forall\,r,s,p\in[m+1,n]\cap\mathbb{N}$; \vskip2mm (A\,2).\,\,
$(I_uI_s)I_p=I_u(I_sI_p)$ \,\, $\forall\,u\in[1,m]\cap\mathbb{N}$\,\,
 $\forall\,s,p\in[m+1,n]\cap\mathbb{N}$.
\vskip2mm

Obviously, the first $m$ basic vectors $\{I_u\}_{u=1}^m$
are idempotents and
form a semi-simple subalgebra of the algebra $\mathbb{A}_n^m$. The
vectors $\{I_r\}_{r=m+1}^n$ form a nilpotent subalgebra of the
algebra $\mathbb{A}_n^m$. The element $1=\sum_{u=1}^mI_u$ is the
unit of $\mathbb{A}_n^m$.

In the cases where $\mathbb{A}_n^m$ has some specific properties,
the following propositions are true.\vskip2mm

\textbf{Proposition 1 \cite{Shpakivskyi-2014}.} \textit{If there exists the unique
$u_0\in[1,m]\cap\mathbb{N}$ such that $I_{u_0}I_s=I_s$ for all
$s=m+1,\ldots,n$,  then the associativity condition \em (A\,2) \em
is satisfied.}\vskip2mm

Thus, under the conditions of Proposition 1,
the associativity condition (A\,1) is only required.
It means that the nilpotent subalgebra of $\mathbb{A}_n^m$ with
the basis $\{I_r\}_{r=m+1}^n$  can be an arbitrary commutative
associative nilpotent algebra of dimension $n-m$. We note that such
nilpotent algebras are fully described for the dimensions
$1,2,3$ in the paper \cite{Burde_de_Graaf}, and some four-dimensional nilpotent algebras
can be found in the papers \cite{Burde_Fialowski},  \cite{Martin}.\vskip2mm

\textbf{Proposition 2 \cite{Shpakivskyi-2014}.} \textit{If all $u_r$ are different in the
multiplication rule \em 3\em ,
 then $I_sI_p=0$ for all $s,p=m+1,\ldots, n$.}\vskip2mm

Thus, under the conditions of Proposition 2,
the multiplication table of the
nilpotent subalgebra of $\mathbb{A}_n^m$ with the basis
$\{I_r\}_{r=m+1}^n$ consists only of zeros, and all associativity
conditions are satisfied.

The algebra $\mathbb{A}_n^m$ contains $m$ maximal ideals
$$\mathcal{I}_u:=\Biggr\{\sum\limits_{r=1,\,r\neq u}^n\lambda_rI_r:\lambda_r\in
\mathbb{C}\Biggr\}, \quad  u=1,2,\ldots,m,
$$
and their intersection is the radical $$\mathcal{R}:=
\Bigr\{\sum\limits_{r=m+1}^n\lambda_rI_r:\lambda_r\in
\mathbb{C}\Bigr\}.$$

Consider $m$ linear functionals
$f_u:\mathbb{A}_n^m\rightarrow\mathbb{C}$ satisfying the
equalities
$$f_u(I_u)=1,\quad f_u(\omega)=0\quad\forall\,\omega\in\mathcal{I}_u\,,
\quad u=1,2,\ldots,m.
$$
Inasmuch as the kernel of functional $f_u$ is the maximal ideal
$\mathcal{I}_u$, this functional is also continuous and
multiplicative (see \cite[p. 147]{Hil_Filips}).


\section{Monogenic functions}

Let us consider the vectors $e_1=1,e_2,\ldots,e_k$ in $\mathbb{A}_n^m$, where $2\leq k\leq 2n$, and these vectors are linearly independent over the field of real numbers
$\mathbb{R}$ (see \cite{Pl-Pukh-Analele}). It means that the equality
$$\sum\limits_{j=1}^k\alpha_je_j=0,\qquad \alpha_j\in\mathbb{R},$$
holds if and only if $\alpha_j=0$ for all $j=1,2,\ldots,k$.

Let the vectors $e_1=1,e_2,\ldots,e_k$ have the following decompositions with respect to
the basis $\{I_r\}_{r=1}^n$:
\begin{equation}\label{e_1_e_2_e_3-k}
e_1=\sum\limits_{r=1}^mI_r\,, 
\quad e_j=\sum\limits_{r=1}^na_{jr}\,I_r\,,\quad a_{jr}\in\mathbb{C},\quad j=2,3,\ldots,k.
\end{equation}

Let $\zeta:=\sum\limits_{j=1}^kx_j\,e_j$, where $x_j\in\mathbb{R}$. It is
obvious that
 $$\xi_u:=f_u(\zeta)=x_1+\sum\limits_{j=2}^kx_j\,a_{ju},\quad u=1,2,\ldots,m.$$
 Let
 $E_k:=\{\zeta=\sum\limits_{j=1}^kx_je_j:\,\, x_j\in\mathbb{R}\}$ be the
linear span of vectors $e_1=1,e_2,\ldots,e_k$ over the field
$\mathbb{R}$.

Let $\Omega$ be a domain in $E_k$. With a domain $\Omega\subset E_k$ we associate the domain  $$\Omega_{\mathbb{R}}:=\Big\{(x_1,x_2,\ldots,x_k)\in\mathbb{R}^k:\,\zeta=\sum\limits_{j=1}^kx_j\,e_j\in\Omega\Big\}$$ in $\mathbb{R}^k$.

We say that a continuous function
$\Phi:\Omega\rightarrow\mathbb{A}_n^m$ is \textit{monogenic}
in $\Omega$ if $\Phi$ is differentiable in the sense of
Gateaux in every point of $\Omega$, i.~e. if  for every
$\zeta\in\Omega$ there exists an element
$\Phi'(\zeta)\in\mathbb{A}_n^m$ such that
\begin{equation}\label{monogennaOZNA}\medskip
\lim\limits_{\varepsilon\rightarrow 0+0}
\left(\Phi(\zeta+\varepsilon
h)-\Phi(\zeta)\right)\varepsilon^{-1}= h\Phi'(\zeta)\quad\forall\,
h\in E_k.\medskip
\end{equation}
$\Phi'(\zeta)$ is the \textit{Gateaux derivative} of the function
$\Phi$ in the point $\zeta$.

Consider the decomposition of a function
$\Phi:\Omega\rightarrow\mathbb{A}_n^m$ with respect to the
basis $\{I_r\}_{r=1}^n$:
\begin{equation}\label{rozklad-Phi-v-bazysi-k}
\Phi(\zeta)=\sum_{r=1}^n U_r(x_1,x_2,\ldots,x_k)\,I_r\,.
 \end{equation}

In the case where the functions $U_r:\Omega_{\mathbb{R}}\rightarrow\mathbb{C}$ are
$\mathbb{R}$-differentiable in $\Omega_{\mathbb{R}}$, i.~e. for every $(x_1,x_2,\ldots,x_k)\in\Omega_{\mathbb{R}}$,
$$U_r\left(x_1+\Delta x_1,x_2+\Delta x_2,\ldots,x_k+\Delta x_k\right)-U_r(x_1,x_2,\ldots,x_k)=
$$
$$=\sum\limits_{j=1}^k\frac{\partial U_r}{\partial x_j}\,\Delta x_j+
\,o\left(\sqrt{\sum\limits_{j=1}^k(\Delta x_j)^2}\,\right), \qquad \sum\limits_{j=1}^k(\Delta x_j)^2\to 0\,,$$
the function $\Phi$ is monogenic in the domain $\Omega$ if
and only if the following Cauchy~-- Riemann conditions are
satisfied in $\Omega$:
\begin{equation}\label{Umovy_K-R-k}
\frac{\partial \Phi}{\partial x_j}=\frac{\partial \Phi}{\partial
x_1}\,e_j \qquad \text{for all}\quad j=2,3,\ldots,k.
\end{equation}

\section{An expansion of the resolvent}

Let $b:=\sum\limits_{r=1}^nb_r\,I_r\in\mathbb{A}_n^m$, where $b_r\in\mathbb{C}$, and we note that $f_u(b)=b_u$, \, $u=1,2,\ldots,m$.
It follows form the Lemmas 1, 3 of \cite{Shpakivskyi-2014} that
\begin{equation}\label{obern-el}
b^{-1}=\sum\limits_{u=1}^m\frac{1}{b_u}\,I_u+
 \sum\limits_{s=m+1}^{n}\sum\limits_{r=2}^{s-m+1}\frac{\widetilde{Q}_{r,s}}
 {b_{u_{s}}^r}\,I_{s}\,.
\end{equation}
where
$\widetilde{Q}_{r,s}$ are determined by the following recurrence
relations:
\begin{equation}\label{Q-hv}
\begin{array}{c}
\displaystyle
\widetilde{Q}_{2,s}:=b_s\,,\qquad
\widetilde{Q}_{r,s}=\sum\limits_{q=r+m-2}^{s-1}\widetilde{Q}_{r-1,q}\,\widetilde{B}_{q,\,s}\,,\; \;\;r=3,4,\ldots,s-m+1,\\
\end{array}
\end{equation}
\begin{equation}\label{B-hv}
\widetilde{B}_{q,s}:=\sum\limits_{p=m+1}^{s-1}b_p \Upsilon_{q,s}^p\,,
\;\;p=m+2,m+3,\ldots,n,
\end{equation}
 and the natural numbers $u_s$ are defined in the
rule  3 of the multiplication table of algebra
$\mathbb{A}_n^m$.

In the next lemma we find an expansion of the resolvent $(te_1-\zeta)^{-1}$.

\vskip2mm
\lem\label{lem_1_rezolv_A_n_m-k} \textit{An expansion of the
resolvent is of the form
 \begin{equation}\label{lem-rez}
(te_1-\zeta)^{-1}=\sum\limits_{u=1}^m\frac{1}{t-\xi_u}\,I_u+
 \sum\limits_{s=m+1}^{n}\sum\limits_{r=2}^{s-m+1}\frac{Q_{r,s}}
 {\left(t-\xi_{u_{s}}\right)^r}\,I_{s}\,
 \end{equation}
$$\forall\,t\in\mathbb{C}:\,
t\neq \xi_u,\quad u=1,2,\ldots,m,$$
where the coefficients\ $Q_{r,s}$ are determined by the following recurrence
relations:
\begin{equation}\label{Q}
\begin{array}{c}
\displaystyle
Q_{2,s}=T_s\,,\quad
Q_{r,s}=\sum\limits_{q=r+m-2}^{s-1}Q_{r-1,q}\,B_{q,\,s}\,,\; \;\;r=3,4,\ldots,s-m+1,\\
\end{array}
\end{equation}
with
\begin{equation}\label{B}
T_s:=\sum\limits_{j=2}^kx_ja_{js}\,, \quad
B_{q,s}:=\sum\limits_{p=m+1}^{s-1}T_p \Upsilon_{q,s}^p\,,
\;\;p=m+2,m+3,\ldots,n,
\end{equation}
 and the natural numbers $u_s$ are defined in the
rule \em 3 \em of the multiplication table of algebra}
$\mathbb{A}_n^m$.\vskip2mm

\textbf{Proof.} Taking into account the decomposition
$$te_1-\zeta=\sum\limits_{u=1}^m(t-\xi_u)I_u-\sum\limits_{r=m+1}^n\sum\limits_{j=2}^kx_ja_{js}\,I_r\,,
$$
we conclude that the relation (\ref{lem-rez}) follows directly from the
equality (\ref{obern-el}) in which instead of $b_u$,\, $u=1,2,\ldots,m$ it should be used the
expansion $t-\xi_u$\,, and  instead of $b_s$\,, $s=m+1,m+2,\ldots,n$ it should be used the
expansion $\sum\limits_{j=2}^kx_ja_{js}$.
 The lemma is proved.

It follows from Lemma \ref{lem_1_rezolv_A_n_m-k} that the points
 $(x_1,x_2,\ldots,x_k)\in\mathbb{R}^k$
corresponding to the noninvertible elements
$\zeta=\sum\limits_{j=1}^kx_j\,e_j$
form the set
  \[M_u^{\mathbb{R}}:\quad\left\{
\begin{array}{r}x_1+\sum\limits_{j=2}^kx_j\,{\rm Re}\,a_{ju}=0,\vspace*{3mm} \\
\sum\limits_{j=2}^kx_j\,{\rm Im}\,a_{ju}=0, \\ \medskip
\end{array} \right. \qquad u=1,2,\ldots,m\]
in the $k$-dimensional space $\mathbb{R}^k$.
Also we consider the set
$M_u:=\{\zeta\in E_k: f_u(\zeta)=0\}$ for $u=1,2,\ldots,m$.
It is obvious that the set $M_u^{\mathbb{R}}\subset\mathbb{R}^k$ is congruent with the set  $M_u\subset E_k$.

\vskip3mm

\section{A constructive description of monogenic functions}
\vskip3mm

We say that a domain $\Omega\subset E_k$ is \textit{convex with respect to the set of directions}
$M_u$ if  $\Omega$ contains the segment $\{\zeta_1+\alpha(\zeta_2-\zeta_1):\alpha\in[0,1]\}$ for all $\zeta_1,\zeta_2\in \Omega$ such that $\zeta_2-\zeta_1\in M_u$.

Denote $f_u(E_k):=\{f_u(\zeta) : \zeta\in E_k\}$. In what follows,
we make the following essential assumption:
$f_u(E_k)=\mathbb{C}$ for all\, $u=1,2,\ldots,m$.
 Obviously, it holds if and only if for every fixed $u=1,2, \ldots, m$
at least one of the numbers $a_{2u}$, $a_{3u},\ldots,a_{ku}$ belongs to
$\mathbb{C}\setminus\mathbb{R}$.
\vskip2mm

\lem\label{lem_1_konstruct_opys_A_n_m-k} \textit{Suppose that a
domain $\Omega\subset E_k$ is convex with respect to the set of directions
$M_u$  and $f_u(E_k)=\mathbb{C}$ for all
 $u=1,2,\ldots, m$. Suppose also that a
function $\Phi:\Omega\rightarrow\mathbb{A}_n^m$ is
monogenic in the domain
  $\Omega$. If points
$\zeta_{1},\zeta_{2}\in\Omega$ such that
$\zeta_{2}-\zeta_{1}\in M_u$,
then}
\begin{equation}\label{Fi(dz1')-Fi(dz2')}
\Phi(\zeta_2)-\Phi(\zeta_1)\in\mathcal{I}_u\,.
\end{equation}

\textbf{Proof.} Inasmuch as  $f_u(E_k)=\mathbb{C}$, then there exists an element $e_2^*\in E_k$
such that $f_u(e_2^*)=i$. Consider the lineal span $E^*:=\{\zeta=xe_1^*+ye_2^*+ze_3^*:x,y,z\in\mathbb{R}\}$ of the vectors $e_1^*:=1, e_2^*,e_3^*:=\zeta_2-\zeta_1$ and denote $\Omega^*:=\Omega\cap E^*$.

Now, the relations $(\ref{Fi(dz1')-Fi(dz2')})$ can be proved in such a way as Lemma 2.1 \cite{Pl-Shp1}, in the proof of which one must take $\Omega^*,f_u,\{\alpha e_3^*:\alpha\in\mathbb{R}\}$ instead of $\Omega_\zeta, f,L$, respectively. Lemma \ref{lem_1_konstruct_opys_A_n_m-k} is proved.

Let a domain $\Omega\subset E_k$ be convex with respect to the set of directions
$M_u$\,, $u=1,2,\ldots, m$. By
$D_u$ we denote that domain in $\mathbb{C}$
 onto which the domain $\Omega$ is mapped by the functional  $f_u$.

We introduce the linear operators $A_u$\,, $u=1,2,\ldots,m$, which
assign holomorphic functions $F_u:\,D_u\rightarrow\mathbb{C}$ to
every monogenic function
$\Phi:\Omega\rightarrow\mathbb{A}_n^m$ by the formula
\begin{equation}\label{def_op_A-k}
F_u(\xi_u)=f_u(\Phi(\zeta)),
\end{equation}
where $\xi_u=f_u(\zeta)\equiv x_1+\sum\limits_{j=2}^kx_j\,a_{ju}$ and
$\zeta\in\Omega$. It follows from Lemma
\ref{lem_1_konstruct_opys_A_n_m-k}  that the value $F_u(\xi_u)$ does
not depend on a choice of a point $\zeta$ for which
$f_u(\zeta)=\xi_u$.

Now, similar to proof of Lemma 5 \cite{Shpakivskyi-2014} can be proved the following statement.

\vskip2mm
\lem\label{lem_2_konstruct_opys_A_n_m-k} \textit{Suppose
that a domain $\Omega\subset E_k$ is  convex with respect to the set of directions
$M_u$ and $f_u(E_k)=\mathbb{C}$
for all $u=1,2,\ldots, m$. Suppose also that for any fixed
$u=1,2,\ldots,m$, a function $F_u:D_u\rightarrow \mathbb{C}$ is
holomorphic in a domain $D_u$ and $\Gamma_u$ is a closed Jordan
rectifiable curve in $D_u$ which surrounds the point $\xi_u$ and
contains no points $\xi_q$, $q=1,2,\ldots, m$,\, $q\neq u$. Then
the function
\begin{equation}\label{lem_5-k}
\Psi_u(\zeta):=I_u\int\limits_{\Gamma_u}F_u(t)(te_1-\zeta)^{-1}\,dt
\end{equation}
is monogenic in the domain
 $\Omega$. }

\vskip2mm

\lem\label{Lem-6-osn-const-op-k} \textit{Suppose that a
domain $\Omega\subset E_k$ is convex with respect to the set of directions
$M_u$  and $f_u(E_k)=\mathbb{C}$ for all
 $u=1,2,\ldots, m$. Suppose also that a function $V:\Omega_{\mathbb{R}}\rightarrow\mathbb{C}$
 satisfies the equalities
 \begin{equation}\label{lem-6-1-k}
\frac{\partial V}{\partial x_2}=\frac{\partial V}{\partial x_1}\,
a_{2u}\,,\quad\frac{\partial V}{\partial x_3}=\frac{\partial V}{\partial x_1}\,
a_{3u}\,,\quad\ldots, \quad \frac{\partial V}{\partial x_k}=\frac{\partial
V}{\partial x_1}\, a_{ku}
\end{equation}
in $\Omega_{\mathbb{R}}$. Then $V$ is a holomorphic function of the variable
$\xi_u=f_u(\zeta)=x_1+\sum\limits_{j=2}^kx_j\,a_{ju}$ in the domain $D_u$.}
\vskip2mm

\textbf{Proof.} We first separate the real and the imaginary part
of the expression
\begin{equation}\label{lem-6-2-k}
\xi_u=x_1+\sum\limits_{j=2}^kx_j\,{\rm Re}\,a_{ju}+i\,
\sum\limits_{j=2}^kx_j\,{\rm Im}\,a_{ju}=:\tau_u+i\eta_u
\end{equation}
and note that the equalities (\ref{lem-6-1-k}) yield
\begin{equation}\label{lem-6-3-k}
\frac{\partial V}{\partial \eta_u}\,{\rm Im}\,a_{2u} =i\,\frac{\partial
V}{\partial \tau_u}\,{\rm Im}\,a_{2u}\,, \quad\ldots,\quad \frac{\partial
V}{\partial \eta_u}\,{\rm Im}\,a_{ku}=i\,\frac{\partial V}{\partial
\tau_u}\, {\rm Im}\,a_{ku}\,.
\end{equation}

It follows from the condition $f_u(E_k)=\mathbb{C}$ that
 at least one of the numbers ${\rm Im}\,a_{2u}$\,, ${\rm Im}\,a_{3u}\,,\ldots,{\rm Im}\,b_u$ is not
equal to zero. Therefore, using (\ref{lem-6-3-k}), we get
\begin{equation}\label{lem-6-4-k}
\frac{\partial V}{\partial \eta_u}=i\,\frac{\partial V}{\partial
\tau_u}\,.
\end{equation}

Now we prove that $V(x'_1,x'_2,\ldots,x'_k)=V(x''_1,x''_2,\ldots,x''_k)$ for points
$(x'_1,x'_2,\ldots,x'_k),(x''_1,x''_2,\ldots,x''_k)\in\Omega$ such that the segment that
connects these points is parallel to a straight line $L_u\subset M_u^{\mathbb{R}}$\,.
To this end we use considerations with the proof of Lemma \ref{lem_1_konstruct_opys_A_n_m-k}.
Since $f_u(E_k)=\mathbb{C}$, then there exists an element $e_2^*\in E_k$
such that $f_u(e_2^*)=i$. Consider the lineal span $E^*:=\{\zeta=xe_1^*+ye_2^*+ze_3^*:x,y,z\in\mathbb{R}\}$ of the vectors $e_1^*:=1$, $e_2^*$, $e_3^*:=\zeta'-\zeta''$, where $\zeta':=\sum\limits_{j=1}^kx'_j\,e_j$\,, \, $\zeta'':=\sum\limits_{j=1}^kx''_j\,e_j$\,, and introduce the denotation $\Omega^*:=\Omega\cap E^*$.

Now, the relation $V(x'_1,x'_2,\ldots,x'_k)=V(x''_1,x''_2,\ldots,x''_k)$ can be proved in such a way as Lemma 6 \cite{Shpakivskyi-2014}, in the proof of which one must take $\Omega^*,\{\alpha e_3^*:\alpha\in\mathbb{R}\}$ instead of $\Omega_\zeta\,, L$, respectively. The lemma is proved.

Thus, a function $V:\Omega_{\mathbb{R}}\rightarrow\mathbb{C}$ of the form
$V(x_1,x_2,\ldots,x_k):=F(\xi_u)$,
 where $F(\xi_u)$ is an arbitrary function holomorphic in
the domain $D_u$\,, is a general solution of the system
(\ref{lem-6-1-k}). The lemma is proved.

\vskip2mm
\theor\label{teo_1_konstruct_opys_A_n_m-k} \textit{Suppose that a
domain $\Omega\subset E_k$ is convex with respect to the set of directions
$M_u$  and $f_u(E_k)=\mathbb{C}$ for all
 $u=1,2,\ldots, m$.  Then every
monogenic function $\Phi:\Omega\rightarrow\mathbb{A}_n^m$
can be expressed in the form
 \begin{equation}\label{Teor--1-k}
\Phi(\zeta)=\sum\limits_{u=1}^mI_u\,\frac{1}{2\pi
i}\int\limits_{\Gamma_u} F_u(t)(te_1-\zeta)^{-1}\,dt+
\sum\limits_{s=m+1}^nI_s\,\frac{1}{2\pi i}\int\limits_
{\Gamma_{u_s}}G_s(t)(te_1-\zeta)^{-1}\,dt,
 \end{equation}\vskip1mm
\noindent where $F_u$ and $G_s$ are certain holomorphic functions in the
domains $D_u$ and $D_{u_s}$, respectively, and $\Gamma_q$ is a
closed Jordan rectifiable curve in $D_q$ which surrounds the point
$\xi_q$ and contains no points $\xi_{\ell}$, $\ell,q=1,2,\ldots,
m$,\,$\ell\neq q$.}\vskip2mm

\textbf{Proof.} We set
 \begin{equation}\label{teor__1-k}
 F_u:=A_u\Phi,\;\;\;u=1,2,\ldots,m.
 \end{equation}
Let us show that the values of monogenic function
\begin{equation}\label{teor__2-k}
\Phi_0(\zeta):=\Phi(\zeta)-\sum\limits_{u=1}^mI_u\,\frac{1}{2\pi
i} \int\limits_{\Gamma_u}F_u(t)(te_1-\zeta)^{-1}\,dt
\end{equation}
belong to the radical $\mathcal{R}$,
i.~e. $\Phi_0(\zeta)\in\mathcal{R}$ for all
$\zeta\in\Omega$.
As a consequence of the equality (\ref{lem-rez}), we have the
equality
$$I_u\,\frac{1}{2\pi i}
\int\limits_{\Gamma_u}F_u(t)(te_1-\zeta)^{-1}\,dt=I_u\,\frac{1}{2\pi
i} \int\limits_{\Gamma_u}\frac{F_u(t)}{t-\xi_u}\,dt+$$
$$
 +\frac{1}{2\pi i}\sum\limits_{s=m+1}^{n}
 \sum\limits_{r=2}^{s-m+1}\int\limits_{\Gamma_u}\frac{F_u(t)Q_{r,s}}
 {\left(t-\xi_{u_{s}}\right)^r}\,dt
 \,I_{s}\,I_u\,,
$$
from which we obtain the equality
\begin{equation}\label{teor__3-k}
f_u\left(\sum\limits_{u=1}^mI_u\,\frac{1}{2\pi i}
\int\limits_{\Gamma_u}F_u(t)(te_1-\zeta)^{-1}\,dt\right)=F_u(\xi_u).
\end{equation}
Operating onto the  equality (\ref{teor__2-k}) by the functional
 $f_u$ and taking into account the relations (\ref{def_op_A-k}), (\ref{teor__1-k}), (\ref{teor__3-k}),  we get the equality
$$f_u(\Phi_0(\zeta))=F_u(\xi_u)-F_u(\xi_u)=0
$$
for all $u=1,2,\ldots,m$, i.~e. $\Phi_0(\zeta)\in\mathcal{R}$.

Therefore, the function $\Phi_0$ is of the form
\begin{equation}\label{Fi_0--1-k}
\Phi_{0}(\zeta)=\sum\limits_{s=m+1}^{n} V_{s}(x_1,x_2,\ldots,x_k)\,I_s\,,
\end{equation}
where $V_{s}:\Omega_\mathbb{R}\rightarrow\mathbb{C}$\,, and the Cauchy~--
Riemann conditions (\ref{Umovy_K-R-k}) are satisfied with
$\Phi=\Phi_0$.
Substituting the expressions (\ref{e_1_e_2_e_3-k}), (\ref{Fi_0--1-k}) into the equality
(\ref{Umovy_K-R-k}), we obtain
\begin{equation}\label{teor__5-1-k}
\begin{array}{c}
\displaystyle
\sum\limits_{s=m+1}^{n} \frac{\partial V_{s}}{\partial x_2}\,I_s=
\sum\limits_{s=m+1}^{n} \frac{\partial V_{s}}{\partial x_1}\,I_s
\sum\limits_{r=1}^n a_{2r}\,I_r\,,
\vspace*{1mm}\\ \vdots \vspace*{3mm}\\
\displaystyle
\sum\limits_{s=m+1}^{n} \frac{\partial V_{s}}{\partial x_k}\,I_s=
\sum\limits_{s=m+1}^{n} \frac{\partial V_{s}}{\partial x_1}\,I_s
\sum\limits_{r=1}^n a_{kr}\,I_r\,.\\
\end{array}
\end{equation}
Equating the coefficients of $I_{m+1}$ in these equalities,
 we obtain the following system of equations for determining the function
 $V_{m+1}(x_1,x_2,\ldots,x_k)$:
$$\frac{\partial V_{m+1}}{\partial x_2}=\frac{\partial V_{m+1}}{\partial x_1}\,
a_{2\, u_{m+1}}\,,
\quad\ldots,\quad
\frac{\partial V_{m+1}}{\partial x_k}=\frac{\partial V_{m+1}}{\partial x_1}\,
a_{k\,u_{m+1}}\,.
$$
It follows from Lemma \ref{Lem-6-osn-const-op-k} that
$V_{m+1}(x_1,x_2,\ldots,x_k)\equiv G_{m+1}(\xi_{u_{m+1}})$, where $G_{m+1}$ is
a function holomorphic in the domain $D_{u_{m+1}}$\,. Therefore,
\begin{equation}\label{teor__4}
\Phi_0(\zeta)=G_{m+1}(\xi_{u_{m+1}})\,I_{m+1}+
\sum\limits_{s=m+2}^{n} V_{s}(x_1,x_2,\ldots,x_k)\,I_s\,.
\end{equation}

Due to the expansion (\ref{lem-rez}), we have the representation
\begin{equation}\label{teor__5}
I_{m+1}\,\frac{1}{2\pi
i}\int\limits_{\Gamma_{u_{m+1}}}G_{m+1}(t)(te_1-\zeta)^{-1}\,dt=
G_{m+1}(\xi_{u_{m+1}})\,I_{m+1}+\Psi(\zeta),
\end{equation}
where $\Psi(\zeta)$ is a function with values in the set
$\big\{\sum_{s=m+2}^n\alpha_s\,I_s:\alpha_s\in\mathbb{C}\big\}$.

Now, consider the function
 $$\Phi_1(\zeta):=\Phi_0(\zeta)-I_{m+1}\,\frac{1}{2\pi
i}\int\limits_{\Gamma_{u_{m+1}}}G_{m+1}(t)(te_1-\zeta)^{-1}\,dt.$$
In view of the relations (\ref{teor__4}), (\ref{teor__5}),
$\Phi_1$ can be represented in the form
$$\Phi_{1}(\zeta)=\sum\limits_{s=m+2}^n \widetilde
V_{s}(x_1,x_2,\ldots,x_k)\,I_s\,,$$
where $\widetilde V_{s}:\Omega_\mathbb{R}\rightarrow\mathbb{C}$\,.

Inasmuch as $\Phi_1$ is a monogenic function in $\Omega$,
the functions $\widetilde V_{m+2},\widetilde
V_{m+3},\dots,\widetilde V_{n}$ satisfy the system
\eqref{teor__5-1-k}, where $V_{m+1}\equiv 0$, $V_s=\widetilde V_{s}$
for $s=m+2,m+3,\ldots,n$. Therefore,  similarly to the function
$V_{m+1}(x_1,x_2,\ldots,x_k)\equiv G_{m+1}(\xi_{u_{m+1}})$, the function
$\widetilde V_{m+2}$ satisfies the equations
$$
\frac{\partial \widetilde V_{m+2}}{\partial x_2}=\frac{\partial \widetilde V_{m+2}}
{\partial x_1}\,a_{2\,u_{m+2}}\,,\quad\ldots,\quad
\frac{\partial \widetilde V_{m+2}}{\partial x_k}=\frac{\partial
\widetilde V_{m+2}}{\partial x_1}\,a_{k\,u_{m+2}}
$$
 and is of the
form $\widetilde V_{m+2}(x_1,x_2,\ldots,x_k)\equiv G_{m+2}(\xi_{u_{m+2}})$,
where $G_{m+2}$ is a function holomorphic in the domain
$D_{u_{m+2}}$\,.

In such a way, step by step, considering the functions
$$\Phi_j(\zeta):=\Phi_{j-1}(\zeta)-I_{m+j}\,\frac{1}{2\pi
i}\int\limits_{\Gamma_{u_{m+j}}}G_{m+j}(t)(te_1-\zeta)^{-1}\,dt$$
 for $j=2,3,\dots,n-m-1$, we get the representation
(\ref{Teor--1-k}) of the function $\Phi$. The theorem is proved.

Taking into account the expansion (\ref{lem-rez}), one can rewrite
the equality (\ref{Teor--1-k}) in the following equivalent form:
$$\Phi(\zeta)=\sum\limits_{u=1}^mF_u(\xi_u)I_u+
\sum\limits_{s=m+1}^{n}\sum\limits_{r=2}^{s-m+1}\frac{1}{(r-1)!}\,
Q_{r,s}\,F_{u_s}^{(r-1)}(\xi_{u_s})\,I_{s}+$$
\begin{equation}\label{dopolnenije-1-1-9-k}+\sum\limits_{q=m+1}^nG_q(\xi_{u_q})I_q+
\sum\limits_{q=m+1}^n\sum\limits_{s=m+1}^{n}
\sum\limits_{r=2}^{s-m+1}\frac{1}{(r-1)!}\,Q_{r,s}\,G_q^{(r-1)}(\xi_{u_q})\,I_{q}
\,I_s\,.
\end{equation}\vskip4mm

Thus, the equalities (\ref{Teor--1-k}) and (\ref{dopolnenije-1-1-9-k})
specify methods to construct explicitly any monogenic functions
$\Phi:\Omega\rightarrow \mathbb{A}_n^m$ using $n$
corresponding holomorphic functions of complex variables.

The following statement follows immediately from the equality
(\ref{dopolnenije-1-1-9-k}) in which the right-hand side is a
monogenic function in the domain $\Pi:=\{\zeta\in
E_k:f_u(\zeta)=D_u,\,u=1,2,\ldots,m\}$.\vskip2mm

\theor\label{teo_pro_naslidky-k} {\it Let a domain $\Omega\subset E_k$ is convex with respect to the set of directions $M_u$  and $f_u(E_k)=\mathbb{C}$ for all
 $u=1,2,\ldots, m$. Then every monogenic function
$\Phi:\Omega\rightarrow \mathbb{A}_n^m$ can be continued
to a function monogenic in the domain $\Pi$.}

The next statement is a fundamental consequence of the equality
(\ref{dopolnenije-1-1-9-k}), and it is true for an arbitrary domain
$\Omega$.\vskip2mm

\theor\label{teo_pro_naslidky2-k} {\it Let $f_u(E_k)=\mathbb{C}$ for all
 $u=1,2,\ldots,m$. Then for every monogenic function
$\Phi:\Omega\rightarrow \mathbb{A}_n^m$ in an arbitrary
domain $\Omega$, the Gateaux $r$-th derivatives
$\Phi^{(r)}$ are monogenic functions in $\Omega$ for all\,
$r$.}

The proof is completely analogous to the proof of Theorem 4
\cite{Pl-Shp1}.

 Using the integral expression
(\ref{Teor--1-k}) of monogenic
function $\Phi:\Omega\rightarrow \mathbb{A}_n^m$ in the case
where a domain $\Omega$ is convex with respect to the set of directions $M_u$\,,\;$u=1,2,\ldots,m$, we obtain the following
expression for the Gateaux $r$-th derivative $\Phi^{(r)}$:
$$\Phi^{(r)}(\zeta)=\sum\limits_{u=1}^mI_u\,\frac{r!}{2\pi i}\int\limits_{\Gamma_u}
F_u(t)\Big((te_1-\zeta)^{-1}\Big)^{r+1}\,dt+$$
$$
+\sum\limits_{s=m+1}^nI_s\,\frac{r!}{2\pi i}\int\limits_
{\Gamma_{u_s}}G_s(t)\Big((te_1-\zeta)^{-1}\Big)^{r+1}\,dt\qquad
\forall\;\zeta\in\Omega\,.\medskip
$$

\section{Remarks}

We note that in the cases where the algebra $\mathbb{A}_n^m$ has some
specific properties (for instance, properties described in
Propositions 1 and 2),
it is easy to simplify the form of the equality
(\ref{dopolnenije-1-1-9-k}). \vskip2mm

\textbf{1.} In the case considered in Proposition 1,
the following equalities hold:
$$u_{m+1}=u_{m+2}=\ldots=u_n=:\eta\,.$$

In this case the representation (\ref{dopolnenije-1-1-9-k}) takes
the form
$$\Phi(\zeta)=\sum\limits_{u=1}^mF_u(\xi_u)I_u+
\sum\limits_{s=m+1}^{n}\sum\limits_{r=2}^{s-m+1}\frac{1}{(r-1)!}\,
Q_{r,s}\,F_\eta^{(r-1)}(\xi_\eta)\,I_{s}+$$
\begin{equation}\label{dopolnenije-1-1-k}+\sum\limits_{s=m+1}^nG_s(\xi_\eta)I_s+
\sum\limits_{q=m+1}^n\sum\limits_{s=m+1}^{n}
\sum\limits_{r=2}^{s-m+1}\frac{1}{(r-1)!}\,Q_{r,s}\,G_q^{(r-1)}(\xi_\eta)\,I_{s}
\,I_q\,.
\end{equation}\vskip2mm

The formula (\ref{dopolnenije-1-1-k}) generalizes representations of
monogenic functions in both three-dimensional harmonic algebras
(see \cite{Pl-Shp1,Pl-Pukh,Pukh}) and specific $n$-dimensional
algebras (see \cite{Pl-Shp-Algeria,Pl-Pukh-Analele}) to the case
of algebras more general form and to a variable of more general form.

\vskip2mm \textbf{2.} In the case considered in Proposition 2, the representation (\ref{Teor--1-k}) takes the form
 \begin{equation}\label{dopolnenije-3-k}
\Phi(\zeta)=\sum\limits_{u=1}^mF_u(\xi_u)I_u+\sum\limits_{s=m+1}^nG_s(\xi_{u_s})I_s+
\sum\limits_{s=m+1}^nT_sF_{u_s}^{\,'}(\xi_{u_s})I_s\,.
\end{equation}

The formula (\ref{dopolnenije-3-k}) generalizes representations of
monogenic functions in both a three-dimensional harmonic algebra
with one-dimensional radical (see \cite{Pl-Pukh}) and semi-simple
algebras (see \cite{Pukh,Pl-Pukh-Analele}) to the case of algebras
more general form and to a variable of more general form.

\vskip2mm \textbf{3.} In the case where $n=m$, the algebra
$\mathbb{A}_n^n$ is semi-simple and contains no nilpotent
subalgebra. Then the formulae (\ref{dopolnenije-1-1-k}),
(\ref{dopolnenije-3-k}) take the form
 $$
 \Phi(\zeta)=\sum\limits_{u=1}^nF_u(\xi_u)I_u\,,
$$
because there are no vectors $\{I_k\}_{k=m+1}^n$.
This formula was obtained in the paper \cite{Pl-Pukh-Analele}.

\section{The relations between monogenic functions and partial differential equations}

Consider the following linear partial differential equation with
constant coefficients:
\begin{equation}\label{dopolnenije----1-k}
\mathcal{L}_NU(x_1,x_2,\ldots,x_k):=\sum\limits_{\alpha_1+\alpha_2\ldots+\alpha_k=N}
C_{\alpha_1,\alpha_2,\ldots,\alpha_k}\,
\frac{\partial^N \Phi}{\partial x_1^{\alpha_1}\,\partial x_2^{\alpha_2}\,\ldots\partial x_k^{\alpha_k}}=0,
\end{equation}\vskip2mm

If a function $\Phi(\zeta)$  is $N$-times
differentiable in the sense of Gateaux in every point of
$\Omega$, then
$$\frac{\partial^{\alpha_1+\alpha_2+\ldots+\alpha_k}\Phi}
{\partial x_1^{\alpha_1}\,\partial x_2^{\alpha_2}\ldots\partial x_k^{\alpha_k}}=$$
$$=
e_1^{\alpha_1}\, e_2^{\alpha_2}\ldots e_k^{\alpha_k}\,\Phi^{(\alpha_1+\alpha_2+\ldots+\alpha_k)}(\zeta)=
e_2^{\alpha_2}e_3^{\alpha_3}\ldots e_k^{\alpha_k}
\,\Phi^{(N)}(\zeta).
$$
Therefore, due to the equality
\begin{equation}\label{dopolnenije----2-1-k}
\mathcal{L}_N\Phi(\zeta)=\Phi^{(N)}(\zeta)
\sum\limits_{\alpha_1+\alpha_2+\ldots+\alpha_k=N}C_{\alpha_1,\alpha_2,\ldots,\alpha_k}\,e_2^{\alpha_2}e_3^{\alpha_3}\ldots e_k^{\alpha_k}\,,
\end{equation}
every $N$-times differentiable in the sense of Gateaux in
$\Omega$ function $\Phi$
satisfies the equation $\mathcal{L}_N\Phi(\zeta)=0$ everywhere in
$\Omega$ if and only if
\begin{equation}\label{dopolnenije----2-k}
\sum\limits_{\alpha_1+\alpha_2+\ldots+\alpha_k=N}C_{\alpha_1,\alpha_2,\ldots,\alpha_k}\,e_2^{\alpha_2}e_3^{\alpha_3}\ldots e_k^{\alpha_k}=0.
\end{equation}
Accordingly, if the condition (\ref{dopolnenije----2-k}) is
satisfied, then the real-valued components ${\rm Re}\,U_k(x_1,x_2,\ldots,x_k)$
and ${\rm Im}\,U_k(x_1,x_2,\ldots,x_k)$ of the decomposition
(\ref{rozklad-Phi-v-bazysi-k})
are solutions of the equation (\ref{dopolnenije----1-k}).

In the case where $f_u(E_k)=\mathbb{C}$ for all $u=1,2, \ldots,
m$, it follows from Theorem \ref{teo_pro_naslidky2-k} that the
equality (\ref{dopolnenije----2-1-k}) holds for every monogenic
function $\Phi:\Omega\rightarrow \mathbb{A}_n^m$.

Thus,
to construct solutions of the equation (\ref{dopolnenije----1-k}) in
the form of components of monogenic
functions,
we must to find $k$ linearly independent over the field $\mathbb{R}$
vectors $(\ref{e_1_e_2_e_3-k})$ satisfying the characteristic
equation (\ref{dopolnenije----2-k}) and to verify the condition:
$f_u(E_k)=\mathbb{C}$ for all $u=1,2,\ldots,m$. Then, the formula
(\ref{Teor--1-k}) gives a constructive description of all mentioned
monogenic functions.

In the next theorem, we assign a special class of equations
(\ref{dopolnenije----1-k}) for which $f_u(E_k)=\mathbb{C}$ for all
$u=1,2,\ldots,m$. Let us introduce the polynomial
 \begin{equation}\label{dopolnenije----51-k}
P(b_2,b_3,\ldots,b_k):=\sum\limits_{\alpha_1+\alpha_2+\ldots+\alpha_k=N}C_{\alpha_1,\alpha_2,\ldots,\alpha_k}\,
b_2^{\alpha_2}\,b_3^{\alpha_3}\ldots b_k^{\alpha_k}.
\end{equation}
\vskip1mm

\theor\label{teo_dopolnenije-dlja-uravn-k} {\it Suppose that
there exist linearly independent over the field $\mathbb{R}$
vectors $e_1=1,e_2,\ldots,e_k$ in $\mathbb{A}_n^m$ of the form
$(\ref{e_1_e_2_e_3-k})$ that satisfy the equality
$(\ref{dopolnenije----2-k})$. If $P(b_2,b_3,\ldots,b_k)\neq0$ for all real $b_2,b_3,\ldots,b_k$, then $f_u(E_k)=\mathbb{C}$ for all $u=1,2,\ldots,m$.}\vskip2mm

\textbf{Proof.} Using the multiplication table of
$\mathbb{A}_n^m$, we obtain the equalities
$$e_2^{\alpha_2}=\sum\limits_{u=1}^ma_{2u}^{\alpha_2}\,I_u+\Psi_\mathcal{R}\,,\quad\ldots,\quad
e_k^{\alpha_k}=\sum\limits_{u=1}^ma_{ku}^{\alpha_k}\,I_u+\Theta_\mathcal{R}\,,
$$
where $\Psi_\mathcal{R}\,,\ldots,\Theta_\mathcal{R}\in\mathcal{R}$. Now
the equality
 (\ref{dopolnenije----2-k}) takes the form
\begin{equation}\label{dopolnenije----4-k}
\sum\limits_{\alpha_1+\alpha_2+\ldots+\alpha_k=N}C_{\alpha_1,\alpha_2,\ldots,\alpha_k}\Biggr(
\sum\limits_{u=1}^ma_{2u}^{\alpha_2}\ldots a_{ku}^{\alpha_k}\,I_u+\widetilde{\Psi}_\mathcal{R}\Biggr)=0,
\end{equation}
where $\widetilde{\Psi}_\mathcal{R}\in\mathcal{R}$. Moreover, due
to the assumption that the vectors $e_1,e_2,\ldots,e_k$ of the form
$(\ref{e_1_e_2_e_3-k})$ satisfy the equality
$(\ref{dopolnenije----2-k})$, there exist complex coefficients
$a_{jr}$ for $j=1,2,\ldots,k,$ $r=1,2,\ldots,n$ that satisfy the equality
(\ref{dopolnenije----4-k}).

It follows from the equality (\ref{dopolnenije----4-k}) that
\begin{equation}\label{dopolnenije----5-k}
\sum\limits_{\alpha_1+\alpha_2+\ldots+\alpha_k=N}C_{\alpha_1,\alpha_2,\ldots,\alpha_k}\,a_{2u}^{\alpha_2}\ldots a_{ku}^{\alpha_k}=0,
\qquad u=1,2,\ldots,m.
\end{equation}

Since $P(b_2,b_3,\ldots,b_k)\neq0$ for all $\{b_2,b_3,\ldots,b_k\}\subset\mathbb{R}$, the equalities
(\ref{dopolnenije----5-k}) can be satisfied only if for each
  $u=1,2,\ldots,m$ at least one of the numbers $a_{2u}$, $a_{3u},\ldots,a_{ku}$ belongs to
 $\mathbb{C}\setminus\mathbb{R}$
 that implies the relation $f_u(E_k)=\mathbb{C}$ for all\, $u=1,2,\ldots,m$.
The theorem is proved.

We note that
if $P(b_2,b_3,\ldots,b_k)\neq0$ for all $\{b_2,b_3,\ldots,b_k\}\subset\mathbb{R}$, then
 $C_{N,0,0,\ldots,0}\neq0$ because otherwise
 $P(b_2,b_3,\ldots,b_k)=0$ for $b_2=b_3=\ldots=b_k=0$.

Since the function $P(b_2,b_3,\ldots,b_k)$ is continuous on $\mathbb{R}^k$, the
condition $P(b_2,b_3,\ldots,b_k)\neq0$ means either $P(b_2,b_3,\ldots,b_k)>0$ or $P(b_2,b_3,\ldots,b_k)<0$ for
all real $b_2,b_3,\ldots,b_k$. Therefore, it is obvious that for any
equation (\ref{dopolnenije----1-k}) of elliptic type, the condition
$P(b_2,b_3,\ldots,b_k)\neq0$ is always satisfied for all $\{b_2,b_3,\ldots,b_k\}\subset\mathbb{R}$.
At the same time, there are equations
 (\ref{dopolnenije----1-k}) for which $P(b_2,b_3,\ldots,b_k)>0$ for all $\{b_2,b_3,\ldots,b_k\}\subset\mathbb{R}$, but
which are not elliptic. For example, such is the equation
$$\frac{\partial^3 u}{\partial x_1^3}+\frac{\partial^3 u}{\partial x_1\partial x_2^2}+
\frac{\partial^3 u}{\partial x_1\partial x_3^2}+\frac{\partial^3 u}{\partial x_1\partial x_4^2}=0
$$
considered in $\mathbb{R}^4$.


\end{document}